\documentclass[11pt]{article}
\usepackage{amssymb}
\usepackage{mathrsfs}
\usepackage{latexsym,amsmath,amssymb,amsfonts,epsfig,graphicx,cite,psfrag}
\usepackage{eepic,color,colordvi,amscd}
\usepackage{ebezier}
\usepackage{verbatim}

\newcommand{\qed}{\hfill $\Box $}
\newcommand{\pf}{\noindent {\bf Proof.} }

\newtheorem{theorem}{Theorem}[section]
\newtheorem{lemma}[theorem]{Lemma}

\begin{document}

\title{Almost perfect matchings in $k$-partite $k$-graphs }

\author{Hongliang Lu\footnote{Partially supported by the National Natural
Science Foundation of China under grant No.11471257}\\
School of Mathematics and Statistics\\
Xi'an Jiaotong University\\
Xi'an, Shaanxi 710049, China\\
\medskip \\
Yan Wang\footnote{Partially supported by NSF grant DMS-1600738 through X. Yu} ~and  Xingxing Yu\footnote{Partially supported by NSF grants DMS-1600738 and CNS-1443894}\\
School of Mathematics\\
Georgia Institute of Technology\\
Atlanta, GA 30332}

\date{May 16, 2017}

\maketitle

\date{}

\maketitle

\begin{abstract}

The minimum co-degree threshold  for a perfect matching in a
$k$-graph with $n$ vertices was determined by R\" odl, Ruci\'nski and Szemer\'edi
for the case when $n\equiv 0\pmod k$. Recently,  Han
resolved the remaining cases when $n \not\equiv 0\pmod k$, establishing a conjecture of R\" odl, Ruci\'nski and
Szemer\'edi. In this paper, we determine the minimum co-degree
threshold  for almost perfect matchings in $k$-partite $k$-graphs,  answering a question of R\"odl and Ruci\'nski.

\end{abstract}

\section{Introduction}

 A {\it hypergraph} $H$
consists of a vertex set $V(H)$ and an edge set $E(H)$ whose members
are subsets of $V(H)$. 
A \emph{matching} in  $H$ is a subset of $E(H)$ consisting
of pairwise disjoint edges.  We use  $\nu(H)$  to denote the maximum size of a matching in  $H$.
A subset $I\subseteq V(H)$ is {\it independent} in $H$ if  $H[I]$, the induced
subhypergraph of $H$ on set $I$, contains no edge.

Let $k$ be a positive integer and $[k] := \{1,2,...,k\}$.
For a set $S$, let ${S\choose k}:=\{T\subseteq S: |T|=k\}$.
A hypergraph $H$ is {\it $k$-uniform} if $E(H)\subseteq {V(H)\choose k}$, and a $k$-uniform hypergraph is also
called a {\it $k$-graph}.
Let $H$ be a $k$-graph. 
For $l\in [k-1]$ and $S\in {V(H)\choose l}$, the {\it neighborhood}
of $S$ in $H$ is the set $N_H(S):=\{ U\in {V(H)\choose k-l} : S \cup U \in E(H)\}$.
The {\it degree} of $S$ in $H$ is $d_H(S):=|N_H(S)|$. The  {\it minimum $l$-degree} of $H$, denoted by $\delta_l(H)$,  is
the minimum degree over  all $l$-subsets of $V(H)$. $\delta_{k-1}(H)$ is also called the  minimum \emph{co-degree}  of $H$.

The minimum co-degree threshold  for a perfect matching in an
$n$-vertex $k$-graph for $n\equiv 0\pmod k$ and large $n$ was determined by R\" odl, Ruci\'nski and Szemer\'edi
\cite{RRS09}.
This threshold function is $n/2-k+C$, where $C \in \{3/2, 2, 5/2, 3\}$, depending
on the parity of $n$ and $k$.
When $n \not\equiv 0\pmod k$, Han \cite{Han15}
determined that $\lfloor n/k \rfloor$ is the minimum co-degree threshold 
for an almost perfect matching in an $n$-vertex $k$-graph, by settling a conjecture in \cite{RRS09}.
Note that when $k \geq 3$ there is a gap between the threshold for a perfect matching and the threshold for an almost perfect matching.
There has been a lot of work on $l$-degree conditions for large matchings, and 
we refer the reader to \cite{TZ12,TZ13} and references therein.

Let $k$ be a positive integer.
A hypergraph $H$ is a {\it $k$-partite $k$-graph} if $V(H)$ admits a
partition $V_1,...,V_k$ such that $|e \cap V_i| = 1$ for all $e \in
E(H)$ and $i \in [k]$. The sets $V_1,\ldots, V_k$ are called {\it partition classes}.
Let $H$ be a $k$-partite $k$-graph with partition classes $V_1,...,V_k$.
A set $S\subseteq V(H)$ is  {\it legal} if $|S \cap V_i| \leq 1$ for $i \in [k]$.
The  {\it minimum $l$-degree} of $H$, denoted by $\delta_l(H)$,  is
the minimum degree over  all legal $l$-subsets of $V(H)$. $\delta_{k-1}(H)$ is also called the  minimum \emph{co-degree}  of $H$.


 K\"uhn and Osthus \cite{Kuh06} showed that  the minimum co-degree
 threshold for a perfect matching in a $k$-partite $k$-graph with $n$
 vertices in each part is between $n/2$ and
 $n/2+\sqrt{2n\log n}$.
Aharoni, Georgakopoulos and Spr\" ussel \cite{AGS09} obtained the
following stronger result: Let $H$ be a $k$-partite $k$-graph with
partition classes $V_1,\ldots,V_k$, each of size $n$. If $d_H(f) >
n/2$ for every legal
$(k-1)$-set $f$ contained in $V(H) - V_1$, and if $d_H(g)\geq n/2$
for every legal $(k-1)$-set $g$ contained in $V(H)-V_2$,
then $H$ has a perfect matching.

For almost perfect matchings in $k$-partitie $k$-graphs, K\"uhn and Osthus \cite{Kuh06}  showed
the following: Let $H$ be a $k$-partite $k$-graph with partition
classes $V_1,\ldots,V_k $, each of size $n$. If
$\delta_{k-1}(H)\geq n/k$, then $\nu(H)\ge n-(k-2)$.
R\"odl and Ruc\'inski \cite{Rod09} asked the following question: Is it
true that, for every $k$-partite $k$-graph $H$ with $n$ vertices in
each partition class,  if $\delta_{k-1}(H) \geq n/k$ then $\nu(H)\ge n-1$?
In this paper, we answer this question in the affirmative (for large
$n$) by proving the following result.
\begin{theorem}
\label{main}
Let $k \geq 3$ be an integer and let $H$ be a $k$-partite $k$-graph with $n$ vertices in each partition class.
Suppose $n$ is sufficiently large.
If $\delta_{k-1}(H) \geq n/k$, then $\nu(H)\ge n-1$.
\end{theorem}

The bound $\delta_{k-1}(H) \geq n/k$ in Theorem~\ref{main} is best possible.
To see this, let $V_1,\ldots, V_k$ be pairwise disjoint sets with $|V_i| = n$ for all $i \in [k]$, 
and let $U_i \subseteq V_i$ for $i \in [k]$ such that $|U_i| = \lfloor (n-1)/k \rfloor$.
Let $H_0$ be the $k$-partite $k$-graph with partition classes $V_1, ..., V_k$
such that, for each legal $k$-set $e$,  $e \in E(H_0)$ if and only if $|e \cap (\bigcup_{i = 1}^k U_i)| \geq 1$.
Clearly, $\delta_{k-1}(H_0) \geq \lfloor (n-1)/k \rfloor$ for large enough $n$.
However, if $n \not\equiv 1 \pmod k$, then $\nu(H_0) \leq |U| = k \cdot \lfloor \frac{n-1}{k} \rfloor < n - 1$.
We remark that for $n \equiv 1 \pmod k$, Han, Zang and Zhao \cite{HZZ} later showed 
that $\lfloor n/k \rfloor$ is the right threshold for $\nu(H)\ge n-1$.

Our proof of Theorem~\ref{main} will be carried out in two steps, similar to the
argument in \cite{Han15}, by considering 
whether or not $H$ contains a large independent set.
In Section 2, we prove Theorem~\ref{main} for the 
case when $H$ has a large independent set, see Lemma~\ref{extremal}.    In Section 3, we prove
Theorem~\ref{main} for the case when $H$ does not contain a large
independent set, see Lemma~\ref{non-extremal}. Theorem~\ref{main} is 
an immediate consequence of Lemmas~\ref{extremal} and \ref{non-extremal}.

\section{Hypergraphs with large independent sets}

In this section, we prove Theorem \ref{main} for the case when $H$ contains a large independent set.
We need a result of Pikhurko \cite{Pik08}.

\begin{lemma}\label{Phk}
Let $k \geq  2$ be an integer and  $l\in [k-1]$, let $\alpha>0$ be a constant and $n$ be a
sufficiently large integer.
 Let $H$ be a $k$-partite $k$-graph
with partition classes $V_1,\ldots,V_k$ and assume  $|V_i|=n$ for $i\in
[k]$. If, for every legal $l$-set $S\subseteq  V_1 \cup \cdots \cup V_l$ and
every legal $(k-l)$-set $S'\subseteq  V_{l+1} \cup \cdots \cup V_k$,
\[
\frac{d_H(S)}{n^{k-l}}+\frac{d_{H}(S')}{n^{l}}>1+\alpha,
\]
then $H$ contains a perfect matching.
\end{lemma}

First, we prove a bound on $\nu(H)$ by slightly weakening the lower bound condition on $\delta_{k-1}(H)$.
This lemma is essentially the same as Theorem 11 in \cite{Kuh06}.

\begin{lemma}\label{small-degree}
Let $n>0,k>0,r\ge 0$ be integers with $n\geq kr+(k-1)$. Let $H$ be  a $k$-partite
$k$-graph with partition classes $V_1,\ldots,V_k$ and let
$|V_i| = n$ for $i \in [k]$. If $\delta_{k-1}(H)\geq r$ then
$\nu(H)\ge kr$.
\end{lemma}

\pf Let $M$ be a maximum matching in $H$. If $|M|\ge
kr$ then $\nu(H)\ge kr$. We may thus assume $|M|<kr$.
Hence, for  $i\in [k]$,
\[
|V_i-V(M)| = n - |M| \geq n - kr + 1 \geq (k-1)+1 = k.
\]
So there exist $k$ pairwise vertex-disjoint legal $(k-1)$-subsets of $V(H)$, denoted as $f_1,\ldots f_k$, such that for $i\in [k]$, 
$$f_i\subseteq \bigcup_{j\in [k]-\{i\}}(V_j-V(M)).$$
Since $M$ is maximum,  $N_H(f_i)\subseteq V(M)$ for $i\in [k]$.

Since $H$ is a $k$-partite $k$-graph,
$N_H(f_i)\cap N_H(f_j)=\emptyset$  for all $i,j\in [k]$ with $i\neq j$.
Hence, since $\delta_{k-1}(H)\geq r$,
$|\bigcup_{i\in [k]}N_H(f_i)|\geq kr.$
Since $|M|<kr$, it follows from the Pigeonhole principle that there
exist  $e\in M$ and distinct  $f_i,f_j$ such that $N_H(f_i)\cap
e\neq\emptyset$ and $N_H(f_j)\cap e \neq \emptyset$.
Let  $v_i\in N_H(f_i)\cap e$ and $v_j\in N_H(f_j)\cap e$. Then
$v_i\neq v_j$ (as $i\ne j$) and,
hence, $(M-\{e\})\cup \{f_i\cup \{v_i\},f_j\cup\{v_j\}\}$ is a  matching in
$H$ whose size is larger than $|M|=\nu(H)$, a contradiction. \qed

\medskip

We now prove Theorem~\ref{main} for the case when $H$ has a large independent set. For
convenience, we adopt  the following concept which is a $k$-partite
version of a concept in \cite{Han15}.
For $0 \leq \gamma < 1$, a $k$-partite $k$-graph $H$ with $n$ vertices in
each partition class is \textit{$\gamma$-extremal} if $H$ contains an
independent set $W$
with at least $(1-\gamma)(k-1)n/k = (1 - 1/k - \gamma + \gamma/k)n$ vertices in each partition class.

\begin{lemma}
\label{extremal}
Let $n,k$ be positive integers with $k\ge 3$ and  $n$  sufficiently large. Let
$\gamma,\varepsilon $ be sufficiently small real numbers such that
 $0<\gamma\ll \varepsilon < 1/(100k^3)$, $\gamma (1-1/k-\gamma)^{k-1} < \varepsilon^2$, and $\varepsilon < (1/k-4k \varepsilon)^{k-1}/100$.
Let $H$ be a $k$-partite $k$-graph with  partition classes
$V_1,\ldots,V_k$ such that $|V_i| = n$ for $i \in [k]$ and $H$ is $\gamma$-extremal.
If $\delta_{k-1}(H)\geq n/k$, then $\nu(H)\ge n-1$.
\end{lemma}

\pf Since $H$ is $\gamma$-extremal,
 there exists an independent set $W$ in $H$ such that $|W\cap
 V_i| \geq {\lceil}(1-1/k-\gamma)n {\rceil}$ for $i\in [k]$.
 So we choose $W$ such that  
   $|W\cap
 V_i|={\lceil}(1-1/k-\gamma)n{\rceil}$ for $i\in [k]$, and let  $W_i=W\cap V_i$. 
Then for $i\in [k]$, $|V_i - W_i| = n - {\lceil}(1-1/k-\gamma)n{\rceil}$; so 
$$  (1/k + \gamma)n-1 =n - (1-1/k-\gamma)n -1 \le |V_i - W_i| \leq n - (1-1/k-\gamma)n = (1/k + \gamma)n.$$ 

For $i\in [k]$, let 
\[
A_i := \{x\in V_i-W_i\ : \ e_H(x,W)\geq \left((1-1/k-\gamma)^{k-1}-\varepsilon\right)n^{k-1}\},
\]
where  $e_H(x,W)$ denotes the number of edges of $H$ containing $x$
and contained in $W\cup \{x\}$.
We proceed by proving four claims. First we give a lower bound on $|A_i|$ for $i\in [k]$.

\medskip

\textit{Claim 1. } $|A_i|\geq (1/k-\varepsilon)n$ for $i\in [k]$.

For, suppose this is not true and, without loss of generality, assume $|A_1| < (1/k-\varepsilon)n$.
Then $|V_1-W_1|-|A_1|\ge (1/k+\gamma)n{-1}-(1/k-{\varepsilon})n=(\varepsilon+\gamma)n{-1}$.
By definition of $A_1$, for each  $x\in V_1-W_1-A_1$,
$e_H(x,W) < \left((1-1/k-\gamma)^{k-1}-\varepsilon\right)n^{k-1}.$
Hence,
\begin{align*}
\sum_{x\in V_1-W_1}e_H(x,W)
&\leq |A_1| \lceil (1-1/k-\gamma)n \rceil ^{k-1}+|V_1-W_1-A_1|\left((1-1/k-\gamma)^{k-1}-\varepsilon\right)n^{k-1}\\
&\leq |V_1 - W_1| \lceil (1-1/k-\gamma)n \rceil ^{k-1} - |V_1-W_1-A_1|\varepsilon n^{k-1}\\
&\le  (1/k+\gamma)n\lceil (1-1/k-\gamma)n \rceil ^{k-1}-((\varepsilon+\gamma)n{-1}) \varepsilon n^{k-1}\\
&=(1/k)n \lceil (1-1/k-\gamma) n \rceil ^{k-1} + (\gamma (1-1/k-\gamma)^{k-1} -
  \varepsilon^2-\varepsilon \gamma )n^k + o(n^k)\\
&< (n/k) \lceil (1-1/k-\gamma) n \rceil ^{k-1},
\end{align*}
where the last inequality holds since $\gamma (1-1/k-\gamma)^{k-1} <
\varepsilon^2$ and $n$ is large.

On the other hand, since $\delta_{k-1}(H)\geq n/k$ and $W$ is independent,
\begin{align*}
\sum_{x\in V_1-W_1}e_H(x,W)
&= \sum_{ \text{legal } (k-1)\text{-set } S\subseteq  \bigcup_{i=2}^kW_i}d_H(S) \\
&\geq \sum_{ \text{legal } (k-1)\text{-set } S\subseteq \bigcup_{i=2}^kW_i}\delta_{k-1}(H) \\
&\geq (n/k) \lceil (1-1/k-\gamma) n \rceil ^{k-1},
\end{align*}
a contradiction. \ensuremath{\blacksquare}

\medskip

Therefore, for $i\in [k]$, we {\color{black}can} choose $A_i'\subseteq A_i$ such that
$|A_i'|=\lceil(1/k-\varepsilon)n\rceil$. Let $A'=\bigcup_{i=1}^k A_i'$.
Then $H-A'$ is a $k$-partite $k$-graph with
$n-\lceil(1/k-\varepsilon)n\rceil$ vertices in each partition
class.
Let $$r=\lceil n/k \rceil-\lceil(1/k-\varepsilon)n\rceil.$$
Then $$1\leq r \leq (n/k + 1) - (1/k-\varepsilon)n  = \varepsilon n +
1.$$
Since $\delta_{k-1}(H) \geq n/k$ and
$|A_i'|=\lceil(1/k-\varepsilon)n\rceil$,  $$\delta_{k-1}(H-A')\geq \lceil n/k \rceil-\lceil(1/k-\varepsilon)n\rceil =r.$$

\medskip

\textit{Claim 2. } There exists a matching $M_0'$ in $H-A'$ such that $|M_0'| \leq k(\varepsilon n + 1)$ and $n' := |V_1| - |M_0'| = k|A_1'|$.

By Lemma \ref{small-degree}, $H-A'$ contains a matching of size at
least $kr$, say  $M_0$. Write $n=mk+s$ with $1\le s\le k$.
Choose $M_0'\subseteq M_0$ such that
  $|M_0'|=k(r-1)+s$.
Then $|M_0'| \leq k(r-1)+k = kr \leq k(\varepsilon n +1)$, and
\begin{align*}
n'=n-|M_0'|=n-k(r-1)-s
=n-k\left(\lceil n/k\rceil-\lceil(1/k-\varepsilon)n\rceil\right)-s+k
=k\lceil(1/k-\varepsilon)n\rceil.
\end{align*}
Hence, $n'=k|A_1'|$. 
\ensuremath{\blacksquare}

\medskip

Let $H' := H-V(M_0')$.
Then $H'$ is a $k$-partite $k$-graph with $n'$ vertices in each partition class.
Let $D_i :=V_i-A_i'-W_i-V(M_0')$ for $i \in [k]$, and let $D = \cup_{i = 1}^k D_i$. Note that for $i\in [k]$, 
$$|D_i|\le |V_i-W_i|-|A_i'|\le (1/k+\gamma)n-(1/k-\varepsilon)n= (\varepsilon+\gamma)n.$$

\medskip

\textit{Claim 3. } There exists a matching $M_1$ in $H'$ such that $D \subseteq V(M_1)$, each edge of $M_1$ contains at least one vertex from $D$, and each edge of $M_1$ contains at most one vertex from $A'$. In particular, $|M_1| \leq |D| \leq k(\varepsilon+\gamma)n$.

We construct such a matching greedily, by starting with the empty matching $F_0 = \emptyset$.
Suppose for some $i \geq 0$, we have constructed a matching $F_i$ such that each edge of $F_i$ contains at least one vertex from $D$ and at most one vertex from $A'$.

We may assume $D \not \subseteq V(F_i)$; for, otherwise, $F_i$ gives
the desired matching $M_1$ for Claim 3. Let $v \in D - V(F_i)$. Then
$v \in D_j$ for some $j\in [k]$ and, without loss of generality,
assume $j \in [k - 1]$.

Since $|F_i| < |D| \leq k(\varepsilon + \gamma)n$ (by construction) and
$|V_l-V(M_0')|=k|A_l'|$ (by Claim 2) for $l\in [k]$, we have, for $l\in [k]$,
$$|V_l - V(M_0') - A_l' - V(F_i)| = k|A_l'| - |A_l'| - |F_i| = (k-1)\lceil(1/k-\varepsilon)n\rceil  - k(\varepsilon+\gamma)n > 0.$$
Thus, let $v_l \in V_l - V(M_0') - A_l' - V(F_i)$ for $l\in [k]-\{j,j+1\}$.

Since $\delta_{k-1}(H) \geq n/k$ and $|M_0'| \leq k(\varepsilon n + 1)$ (by Claim 2), we have
$$d_{H' - V(F_i)}(\{v\} \cup  \{v_l: l \in [k] - \{j, j+1\}\} ) \geq n/k - k(\varepsilon n + 1) - k(\varepsilon+\gamma)n > n/(2k) > 0 $$
as $\varepsilon <1/(100k^3)$ and $n$ is sufficiently large. 
Thus there exists a vertex $v_{j+1} \in V_{j+1} - V(M_0') - V(F_i)$ such that $e_{i+1} := \{v\} \cup \{v_l: l \in [k] - \{j\}\} \in E(H)$.
Clearly, $\{v\} \subseteq e_{i+1}\cap D$ and $e_{i+1}\cap A'\subseteq\{v_{j+1}\}$. 

Let $F_{i+1} := F_i \cup \{e_{i+1}\}$.
Then, each edge of $F_{i+1}$ uses at least one vertex from $D$ and at most one vertex from $A'$. 
We may continue this process until we obtain a matching $F_t$ in $H'$ such that $D\subseteq V(F_t)$, each edge of $F_t$ 
uses at least one vertex from $D$,  and each edge of $F_t$ uses at most one vertex from $A'$. 

Clearly, $t \leq k(\varepsilon+\gamma)n$ as $|D| \leq k(\varepsilon+\gamma)n$.
So $M_1 := F_t$ is the desired matching. \ensuremath{\blacksquare}

\medskip

Let $W_i' := W_i - V(M_0') - V(M_1)$ for $i \in [k]$. 

\medskip

\textit{Claim 4. } There exist  $l\in [k-1]$  and matching $M_2$ in $H' - V(M_1)$ such that $|M_2|\le |M_1|$ and
$$ (k-1) \left(\sum_{i=1}^k |A_i'-V(M_1)-V(M_2)|-l\right)= \left(\sum_{i=1}^k |W_i'-V(M_2)|\right)-(k-l).$$



Since $|V_i - V(M_0')| = k |A_i'|$ (by Claim 2), $|V_i - V(M_0') - A_i'| = (k-1)|A_i'|$.
So $(k-1) \sum_{i = 1}^k |A_i'| = \sum_{i = 1}^k |V_i - V(M_0') -
A_i'|$.  
Moreover, since each edge of $M_1$ contains
at most one vertex of $A'=\cup_{i=1}^k A_i'$,  $(k-1) \sum_{i=1}^k |A_i' \cap V(M_1)| =  (k-1) |A' \cap V(M_1)| \leq |(V(H) - V(M_0') - A') \cap V(M_1)| = \sum_{i=1}^k |(V_i - V(M_0') - A_i') \cap V(M_1)|$.
Thus, 
\begin{equation*}\label{sum-M2-up}
\begin{split}
(k-1)\sum_{i=1}^k |A_i'-V(M_1)|  &= (k-1) \sum_{i=1}^k (|A_i'| - |A_i' \cap V(M_1)|) \\
&= (k-1) \sum_{i=1}^k |A_i'| -(k-1)\sum_{i=1}^k|A_i'\cap V(M_1)|\\
&\geq \left(\sum_{i=1}^k|V_i-V(M_0')-A_i'|\right) - \left( \sum_{i=1}^k |(V_i - V(M_0') - A_i') \cap V(M_1)| \right)\\
&= \sum_{i=1}^k|V_i-V(M_0')-A_i'-V(M_1)|.
\end{split}
\end{equation*}

Recall that $D_i\subseteq V(M_1)$; so we have 
$W_i'= V_i-V(M_0')-A_i'-V(M_1).$
Therefore,  we get 
\begin{align*}
(k-1)\sum_{i=1}^k |A_i'-V(M_1)|\ge \sum_{i=1}^k |W_i'|
\end{align*}

We now construct the desired matching $M_2$ greedily, starting with the empty matching $R_0 = \emptyset$. 
By the above expression, we have  
\begin{align*}
\Delta_0:= (k-1) \sum_{i=1}^k |A_i'-V(M_1)-V(R_0)|-\sum_{i=1}^k |W_i'-V(R_0)| \ge 0.
\end{align*}
Moreover, since $(k-1) \sum_{i = 1}^k |A_i'| = \sum_{i = 1}^k |V_i - V(M_0') -
A_i'|$, 
$$\Delta_0 \le \sum_{i=1}^k |V(M_1)\cap V_i| \le k|M_1|.$$
Suppose we have constructed a matching $R_j$  in $H'-V(M_1)$, for some integer $j\ge 0$, such that 
\begin{align*} 
\Delta_j:=(k-1) \sum_{i=1}^k |A_i'-V(M_1)-V(R_j)| -\sum_{i=1}^k |W_i'-V(R_j)| \ge 0 
\end{align*}

If for some $l\in [k-1]$
\begin{align*} 
(k-1) \left(\sum_{i=1}^k |A_i'-V(M_1)-V(R_j)|-l\right) = \left(\sum_{i=1}^k |W_i'-V(R_j)|\right)-(k-l),
\end{align*}
then $M_2 := R_j$ gives the desired matching for Claim 4. 

So we may assume that  for each $l\in [k-1]$,
\begin{align*} 
f(l) := 
(k-1) \sum_{i=1}^k |A_i'-V(M_1)-V(R_j)| - \sum_{i=1}^k |W_i'-V(R_j)| -k(l-1) =\Delta_l-k(l-1)\ne 0.
\end{align*}
Note that $f(l) \equiv 0 \pmod k$, since $f(l)$ can be written as
$$k (\sum_{i=1}^k |A_i'-V(M_1)-V(R_j)| ) - \sum_{i=1}^k (|A_i'-V(M_1)-V(R_j)| + |W_i'-V(R_j)|) - k(l-1),$$
and $|A_i'-V(M_1)-V(R_j)| + |W_i'-V(R_j)|$ are the same for all  $i \in [k]$.
Moreover, $f(l+1)=f(l)-k$ for all integers $l$. 

We claim that $f(k)\ge 0$. To see this, we first note that
$f(1)=\Delta_j> 0$ (as $\Delta_j\ge 0$ and $\Delta_j\ne 0$). Hence, since $f(1)\equiv 0 \pmod k$, 
$f(1)\ge k$. This implies $f(2)=f(1)-k\ge 0$; so $f(2)\ge k$ as $f(2)\ne 0$. Continuing this process and using the assumption 
that $f(l)\ne 0$ for $l\in [k-1]$, we conclude that $f(k)\ge 0$.

Let $v_s \in A_s' - V(M_1) - V(R_j)$ for $s \in [k - 1]$. 
Then by Claims 2 and 3 and by the assumption that $|R_j|\le |M_1|$, we have 
\begin{equation*}
\begin{split}
d_{H' - V(M_1) - V(R_j)}(\{v_1, ..., v_{k-1}\})
&\geq \delta_{k-1}(H' - V(M_1) - V(R_j)) \\
&\geq n/k - |M_0'| - |M_1| - |R_j| \\
&\geq n/k - k(\varepsilon n + 1) - k(\varepsilon + \gamma)n - k(\varepsilon + \gamma)n \\
&> 1.
\end{split}
\end{equation*}
Hence, there exists $v_k \in V_k - V(M_0') - V(M_1) - V(R_j)$ such that $e_{j+1} := \{v_1,v_2,...,v_k\} \in E(H' - V(M_1) - V(R_j))$.
Let $R_{j+1} := R_j \cup \{e_{j+1}\}$, which is a matching in
$H'-V(M_1)$. Note that 
\begin{align*} 
\Delta_{j+1}:=(k-1) \sum_{i=1}^k |A_i'-V(M_1)-V(R_{j+1})|-\sum_{i=1}^k
  |W_i'-V(R_{j+1})| \ge f(k)\ge 0.
\end{align*}

We claim that $|R_{j+1}|\le |M_1|$. 
Note from our construction, $|e_{j+1}\cap A_s'|=1$ for $s\in [k-1]$ and $|e_{j+1}\cap W_k'|\le 1$. So 
$$\Delta_{j+1} \leq \Delta_j - (k-1)(k-1) + 1 = \Delta_j - k(k -2) \leq \Delta_j - k,$$
as $k\ge 3$. Since $\Delta_{j+1}\ge 0$ and $\Delta_0\le k|M_1|$, we see that $|R_{j+1}|\le |M_1|$. 

Clearly, by continuing this process, we obtain a matching $F_t$ which gives the desired matching for Claim 4.  \ensuremath{\blacksquare}

\medskip

Let $L\subseteq V(H')-V(M_1\cup M_2)$ be a legal $k$-set such that $|L \cap ( A' - V(M_0') - V(M_1) - V(M_2) )| = l$.
For $i\in [k]$, let $W_i'' := W_i' - V(M_2) -L$ and $A_i'' := A_i'- V(M_0') - V(M_1) - V(M_2) -L$. 
Since $|L \cap A_i| = 1$, 
\begin{align*}\label{size-A1''}
|A_i''| 
&\geq |A_i'| - |M_0'| - |M_1| - |M_2| - 1 \\ 
&\geq \lceil (1/k - \varepsilon)n \rceil - k(\varepsilon n  + 1) - k(\varepsilon + \gamma)n - k(\varepsilon + \gamma)n - 1 \\
&> (1/k - 4k \varepsilon)n.
\end{align*}

Let  $n'' := |W_1''| + |A_1''|$, and $H'' :=H'-L- V(M_1) - V(M_2)$.
Then by Claim 4, $(k-1)\sum_{i=1}^k |A_i''|=\sum_{i=1}^k |W_i''|$, which gives 
$k \sum_{i=1}^k |A_i''| = \sum_{i=1}^k (|W_1''| + |A_1''|) = kn''$; so 
$$\sum_{i=1}^k |A_i''|=n''.$$ 
Hence, for $i \in [k]$, 
$$|W_i''| = n'' - |A_i''| = \sum_{j \in [k] - \{i\}} |A_j''|.$$
Therefore, 
$V(H'')$ admits a partition $F_1,\ldots,F_k$, such that
$F_i\cap V_i=A_i''$ for $i\in [k]$ and, for any $j\in [k]-\{i\}$,
$F_i\cap V_j\subseteq  W_i''$ and
$|F_i\cap V_j|=|A_i''|$. Note that $H[F_i]$ is a $k$-partite $k$-graph with partition classes
$F_1\cap V_1, \ldots, F_1\cap V_k$. 

\medskip

If, for all $i\in [k]$, $H[F_i]$ has a perfect matching, say $M_i'$, then
$M_0' \cup M_1 \cup M_2 \cup (\bigcup_{i=0}^k M_i')$ is a perfect matching in $H-L$, showing 
that $\nu(H)\ge n-1$. 

Thus,  it suffices to show, without loss of generality, that $H[F_1]$
contains a perfect matching.  We  use Lemma~\ref{Phk}; so we need to  bound
$d_{H[F_1]}(x)$ and $d_{H[F_1]}(S)$ for $x\in A_1''$ and legal $(k-1)$-sets
$S\subseteq \bigcup_{j=2}^k(F_1\cap V_j) $.

\medskip

Let $x\in A_1''$. Since $A_1'' \subseteq A_1$, the number of legal
$(k-1)$-subsets of $W$ not forming an edge with $x$ is at most
$$({\lceil}(1-1/k-\gamma)n{\rceil})^{k-1}-((1-1/k-\gamma)^{k-1}-\varepsilon)n^{k-1}{<1.1}\varepsilon n^{k-1}.$$
Note that the number of legal $(k-1)$-subsets of $\bigcup_{j=2}^k(F_1\cap V_j)$ 
containing at least one vertex  from $L \cup V(M_0') \cup V(M_1) \cup V(M_2)$ is at most
$$(|L| + |V(M_0')| + |V(M_1)| + |V(M_2)| )|A_1''|^{k-2} \leq 4k(\varepsilon + \gamma)n|A_1''|^{k-2}.$$
Hence,
$$d_{F_1}(x)
\geq |A_1''|^{k-1} - 4k(\varepsilon + \gamma)n|A_1''|^{k-2} - {1.1}\varepsilon n^{k-1} \geq 0.9 |A_1''|^{k-1},$$
where the last inequality holds because 
$ \lceil (1/k - \varepsilon)n \rceil > |A_1''|> (1/k - 4k \varepsilon)n$, 
 $0<\gamma <
\varepsilon<1/(100k^3)$, and $\varepsilon < \frac{1}{100} (1/k-4k \varepsilon)^{k-1}$.

On the other hand, for any legal $(k-1)$-set $S\subseteq \bigcup_{j=2}^k(F_1\cap V_j) $ and
large $n$,
we have
$$d_{F_1}(S) \geq n/k- |M_0'| - |M_1| - |M_2| - 1 \geq n/k- k(\varepsilon n + 1)- 2k(\varepsilon+\gamma )n - 1
> (1-4k^2\varepsilon)(n/k).$$

So for any $x\in A_1''$ and any legal $(k-1)$-set $S\subseteq  \bigcup_{j=2}^k(F_1\cap V_j) $,
\[
\frac{d_{F_1}(x)}{({\color{black}|A_1''|)}^{k-1}}+\frac{d_{F_1}(S)}{|A_1''|}
\geq 0.9+\frac{d_{F_1}(S)}{n/k}>0.9+(1-4k^2\varepsilon)>3/2,
\]
since $\varepsilon <1/(100k^3)$. By Lemma \ref{Phk}, $F_1$ contains a perfect matching $M_1$.
 \qed

\section{Hypergraphs without large independent sets}

In this section, we prove Theorem \ref{main} for the case when $H$
does not contain large independent sets.
First, we consider the case when
the minimum co-degree is slightly below $n/k$.

\begin{lemma}
\label{almost-matching-non-extremal-new}
Let $\beta$ be a constant with $0 < \beta < 1/2k^2$ and let $k>0$ be an integer.
Let $H$ be a  $k$-partite $k$-graph with partition classes  $V_1,\ldots,V_k$ such that $|V_i| = n$ for $i \in [k]$.
Suppose $\delta_{k-1}(H) > (1-\beta) (n/k)$ and $H$ is not $(k^2\beta)$-extremal.
Then $\nu(H)\ge n-k^2$.
\end{lemma}

\pf Let $M:=\{e_1,e_2,...,e_m\}$ be a maximum matching in $H$; so
$m=\nu(H)$.
Let $U := V(H) - V(M)$. We may assume $m<n-k^2$. Then  $|U\cap V_i|
> k^2$ for $i\in [k]$. By the maximality of
$|M|$,  $U$ is independent in $H$.
Thus, there exist $k^2$ pairwise disjoint legal sets $A_i^j\in {U\choose
  k-1}$, $i,j \in [ k]$, such that $A_i^j \cap V_j =
\emptyset$.
For $j \in [k]$, let  $$C_j:=\{v \in V_j : A_i^j \cup \{v\} \in E(H)
\mbox{ for some } i\in [k]\}$$
and $$D_j: =\{v \in V_j : A_i^j
\cup \{v\} \in E(H) \mbox{ for all $i\in [k]$}\}.$$
 Clearly,  $D_j \subseteq C_j$ for $j \in [k]$.
Let $C= \bigcup_{j \in [k]} C_j$ and $D= \bigcup_{j \in [k]} D_j$.

We claim that for $i\in [m]$, $|C\cap e_i|\leq 1$.
For, otherwise, suppose  $|C \cap e_l| \geq 2$ for some $l\in [m]$.
Let $x,y \in C\cap e_l$ be distinct.
By definition of $C$, there exist sets $A_{i_x}^{j_x}$ and
$A_{i_y}^{j_y}$ such that $A_{i_x}^{j_x} \cup \{x\} \in E(H)$ and
$A_{i_y}^{j_y} \cup \{y\} \in E(H)$. Since $x,y\in e_l$ are distinct,
$A_{i_x}^{j_x}\ne A_{i_y}^{j_y}$; so $A_{i_x}^{j_x}\cap A_{i_y}^{j_y}=\emptyset$.
Hence, $(M -\{ e_l\}) \cup \{A_{i_x}^{j_x} \cup \{x\}, A_{i_y}^{j_y}
\cup \{y\} \}$ is a  matching in $H$, whose size is larger than
$|M|=\nu(H)$, a contradiction.

Since $\delta_{k-1}(H) > (1-\beta) (n/k)$, $|C_j| > (1-\beta)(n/k)$
for $j \in [k]$. Hence, for any $l\in [k]$, $$\sum_{j=1}^k|C_j|>(k-1)(1-\beta)(n/k)+|C_l|.$$
On the other hand, since $|C\cap e_j|\le 1$ and $C_j\cap U=\emptyset$
for $j\in [m]$ (as $U$ is independent), $$\sum_{j=1}^k |C_j|\leq m< n - k^2 <
n.$$
Therefore, for $l\in [k]$,
$$|C_l|<n-(k-1)(1-\beta)(n/k)=(1+(k-1)\beta )(n/k).$$

Now we derive  a lower bound on $|D_j|$ for $j \in [k]$.
Since $\delta_{k-1}(H) > (1-\beta) (n/k)$, 
\begin{align*}
\sum_{i=1}^k d_H(A_i^j)> k(1-\beta)(n/k)
\end{align*}
for $j \in [k]$.
On the other hand, for $j\in [k]$, $$\sum_{i=1}^k d_H(A_i^j)\leq
k|D_j|+(k-1)(|C_j|-|D_j|)=|D_j|+(k-1)|C_j|.$$
Since $|C_j|<(1+(k-1)\beta )(n/k)$, we have for $j\in [k]$,  $$\sum_{i=1}^k d_H(A_i^j)< |D_j|+(k-1)(1+(k-1)\beta )(n/k).$$
Thus, for $j\in [k]$,
\[
|D_j|> k(1-\beta)(n/k)-(k-1)(1+(k-1)\beta )(n/k)> (1-k^2\beta)(n/k).
\]

Let $V_D := \bigcup_{i\in [m],e_i \cap D \neq \emptyset}e_i$.
Then, since $|D \cap e_i| \leq 1$ for $i \in [m]$,
$$|(V_D \cap V_j) - D_j| \geq \sum_{i \in [k] - \{j\}} |D_i| > (k-1)(1-k^2\beta) (n/k)=(1-k^2\beta) (k-1)(n/k)$$
 for $j \in [k]$.
Since $H$ is not $(k^2\beta)$-extremal (recall this definition from
the paragraph preceding Lemma~\ref{extremal}, $H[V_D - D]$ contains at least one edge, say $e_0$.
Let $e_{m_1}, ...,e_{m_l}$ be the edges  in $M$ intersecting
$e_0$, where  $l\in [k]$. For each $s\in [l]$, since $e_{m_s}\cap e_0\ne \emptyset$
and $e_0\subseteq V_D-D$, $e_{m_s}\cap e_0\subseteq V_D-D$, which, together with the definition of $V_D$, implies $e_{m_s}\cap D\ne \emptyset$. So
let $v_{m_s}\in e_{m_s} \cap D$ for  $s \in [l]$.

Since $l\le k$,
$|\{v_{m_1},\ldots, v_{m_l}\}\cap V_j|\le k$  for $j\in [k]$. Hence, it follows from the  definition of $D$ that
there exist pairwise distinct (hence, disjoint) $A_{i_1}^{j_1},A_{i_2}^{j_2},...,A_{i_l}^{j_l}$ such that $A_{i_s}^{j_s} \cup \{v_{m_s}\} \in E(H)$ for $s \in [l]$.
Therefore, $(M - \{e_{m_s}:s\in [l]\}) \cup \{A_{i_s}^{j_s} \cup
\{v_{m_s}\}: s\in [l]\} \cup \{e_0\}$ is a matching in $H$ whose size
is larger than $|M|=\nu(H)$, a contradiction.
\qed

\medskip

Next, we introduce a convenient concept which we will use to augment a
matching. Let $H$ be a $k$-partite $k$-graph with  partition classes
$V_1,\ldots,V_k$ such that $|V_i| = n$ for $i \in [k]$. A set $S\in
{V(H)\choose k+1}$ is said to be of \textit{type $j$} if $|S\cap V_j|=2$ and
$|S\cap V_i|=1$ for $i\in [k]-\{j\}$. For a set $S$ of type $j$, an edge $e\in E(H)$ is said to
be {\it $S$-absorbing} if the following holds: there exist $r \in [k]
- \{j\}$ and $v\in S\cap V_j$ such that $e \cap S = \emptyset$, $S_e:=(S-\{v\}-V_r)\cup (e\cap
V_r)\in E(H)$ and $e_S:=(e-V_j-V_r)\cup \{v\}\cup (S\cap V_r)\in
E(H)$.

Note that if $M$ is a matching in $G$ and $e\in M$ is $S$-absorbing
for some $S$ of type $j$ such that $S \cap V(M)  =\emptyset$ then
$(M-\{e\})\cup \{S_e,e_S\}$ is also a matching in $H$.

We now show that minimum co-degree conditions on three types of legal
$(k-1)$-sets will guarantee many $S$-absorbing edges for a set $S$ of type $j$.

\begin{lemma}
\label{absorbing-counting}
Let $c$ be a real number with $0<c<1$, let  $k,n$ be integers with $k\ge 3$ and $n$ sufficiently large, and
let $H$ be a $k$-partite $k$-graph with partition classes
$V_1,\ldots,V_k$ such that $|V_i| = n$ for $i \in [k]$.
Let $r,s,t\in [k]$ be distinct and assume $\delta_{k-1}(e) \geq cn$ for any  legal
$e\subseteq  \bigcup_{j \in [k] - \{i\} } V_j$ for $i \in \{r,s,t\}$.
Then for any  $S$ of type $j$ with $j\in
\{r,s,t\}$, the number of $S$-absorbing edges in $H$  is at least $c^3n^k/2$.
\end{lemma}

\pf
Without loss of generality, we may assume that $j=s=1$, $r=2$ and $t=k$.
Let $S = \{v_1,v_2,\ldots,v_k,{\color{black}v_1'}\}$ such that ${\color{black}v_1'\in V_1}$ and $v_i\in
V_i$ for $i\in [k]$. So $S$ is a set of type $j=1$. We count
$S$-absorbing edges $\{u_1,\ldots, u_k\}$ by choosing vertices
$u_i\in V_i$ for $i\in [k]$.

First, we choose $u_i \in V_i - \{v_i\}$ for $i\in [k]-\{1,2,k\}$
arbitrarily. So the number of choices for $\{u_i: in [k]-\{1,2,k\}$ is $(n-1)^{k-3}$.

Note that  $\delta_{k-1}(e) \geq cn$ for any legal
$e\subseteq \bigcup_{j \in [k] - \{k\} } V_j$. Hence, for a given choice
of  $\{u_i: in [k]-\{1,2,k\}$ is $(n-1)^{k-3}$, we have at least $cn - 1$ choices of $u_k \in V_k - \{v_k\}$ such
that $\{v_1',v_2,u_3,\ldots,u_k\} \in E(H)$.

Since $\delta_{k-1}(e) \geq cn$ for any legal   $e\subseteq \bigcup_{j \in [k] - \{2\} } V_j$,
we have at least  $cn - 1$
choices for $u_2 \in V_2 - \{v_2\}$ such that
$\{v_1,u_2,v_3, \ldots,v_k\} \in E(H)$.

Again, note that $\delta_{k-1}(e) \geq cn$ for any legal $e \subseteq\bigcup_{j \in [k] - \{1\} }
V_j$. So for a given choice of $\{u_2,\ldots, u_k\}$, we have at least $cn - 2$ choices for $u_1 \in V_1 -
\{v_1,v_1'\}$ such that $\{u_1,u_2,\ldots,u_k\} \in E(H)$.

It is easy to verify that each $\{u_1,u_2,u_3,\ldots,u_k\}$ chosen above is $S$-absorbing.
Therefore, the number of $S$-absorbing edges in $H$ is at least $$(cn -
1)^2(cn - 2)(n-1)^{k-3},$$ 
which is at least  $c^3n^k/2$ for large $n$.
\qed

\medskip

We need Chernoff bounds in the lemma below, whose proof can be found in \cite{mitzenmacher}.

\begin{lemma}
\label{chernoff}
Suppose $X_1, ..., X_n$ are independent random variables taking values in $\{0, 1\}$. Let $X$ denote their sum and $\mu = \mathbb{E}[X]$ denote the expected value of $X$. Then for any $0 < \delta \leq 1$,
$$\mathbb{P}[X \leq (1-\delta) \mu] < e^{-\frac{\delta^2 \mu}{2}}$$
\end{lemma}

Analogous to Fact 2.3 in \cite{Rod09},
we prove a stronger version of the absorbing lemma for $k$-partite $k$-graphs.

\begin{lemma}
\label{absorbing}
For any constant $c > 0$, there exists an integer  $n_0 > 0$ with the following
property: If $H$ is a $k$-partite $k$-graph with partition classes
$V_1,\ldots,V_k$ such that $|V_i| = n \geq n_0$ (for $i \in [k]$) and
 $\delta_{k-1}(H) \geq cn$,
then
there is a matching $M'$ in $H$ such that $|M'|\le (32/c^3)(k+2)\log n$ and,
 for any $j\in [k]$ and any set $S$ of type $j$, at least
 $4(k+2)\log n$ edges in $M'$ are $S$-absorbing.
\end{lemma}

\pf Let $C = 32(k+2)/c^3$. Form the set $M'\subseteq E(H)$ by
choosing each edge of $H$  independently and  uniformly at random with probability $p = (C/2) n^{-k} \log n $.
Thus, $\mathbb{E}[|M'|] =|E(H)|p\le n^k p = (C/2) \log n$.

The number of ordered intersecting pairs of edges in $E(H)$ is at most
$n^k (k n^{k-1})=kn^{2k-1}$; so
the expected number of ordered intersecting pairs of edges in $M'$ is at most $$kn^{2k-1} p^2 = kC^2 \log^2 n/(4n) = o(1).$$
By Markov's inequality, with probability strictly larger than $1/3$, $M'$ is a matching of size at most $C \log n$.

For a set $S$ of type $j$ with $j\in [k]$, let $X_S$ denote the number of $S$-absorbing edges in $M'$.
Then by Lemma \ref{absorbing-counting}, we have
$$\mathbb{E}[X_S] \geq (c^3n^k/2) p = 8(k+2)\log n.$$
By Lemma~\ref{chernoff},
$$\mathbb{P}[X_S \leq \ \mathbb{E}[X_S] /2] \leq \exp(-
\mathbb{E}[X_S]/8) = \exp(-(k+2)\log n) = n^{-(k+2)}.$$

Since there are at most $kn^{k+1}$ sets $S$ of type $j$ for all $j\in
[k]$, it follows from union bound that,  with probability strictly larger than $1/4$,
 $X_S\ge \mathbb{E}[X_S]/2 \geq 4(k+2)\log n$  for all sets
 $S$ of type $j$ and for all $j\in [k]$. Thus, the desired $M'$
 exists.
\qed

\medskip

We now prove Theorem~\ref{main} for the case when $H$ does not
have  large independent sets. Recall the definition of
$\gamma$-extremal from the paragraph preceding Lemma~\ref{extremal}.

\begin{lemma}
\label{non-extremal}
Let $\gamma$ be a real number with $ 0 < \gamma < 1/2$, let  $k,n$ be
integers with $k\ge 3$ and $n$ sufficiently large, and let $H$ be a
$k$-partite $k$-graph with partition classes $V_1,\ldots,V_k$ such that $|V_i| = n$ for $i \in [k]$.
Suppose $H$ is not $\gamma$-extremal.
If $\delta_{k-1}(H) \geq n/k$, then $\nu(H)\ge n-1$.
\end{lemma}

\pf
Applying Lemma \ref{absorbing} to $H$ with $c = 1/k$, we obtain a
matching $M'$ in $H$ such that  $|M'|\le 64k^4 \log n$ and,
for any set $S$ of any type $j\in [k]$, $M'$ has at least $4(k+2)\log n$ $S$-absorbing edges.

Let $H':= H - V(M')$ and $n' := |V(H') \cap V_1|$. Then
$n' = n - |M'| \geq n - 64k^4 \log n$ and
$$\delta_{k-1}(H') \geq \delta_{k-1}(H) - |M'| \geq n/k - 64k^4 \log n
\geq (1-\gamma/(2k^2))n'/k,$$
as $n$ is sufficiently large.

Since $(1-\gamma/2)(k-1)n'/k>(1-\gamma)(k-1)n/k$ (as $n$ is large) and
$H$ is not $\gamma$-extremal,    $H'$ is not $\gamma/2$-extremal.
So by applying Lemma \ref{almost-matching-non-extremal-new} to $H'$ with
$\beta = \gamma/(2k^2)$, we obtain a matching $M''$ in $H'$
with $|M''| \geq n'-k^2$.
Let $M_1':=M'$, $M_1 := M_1' \cup M''$ and $U_1 := V(H) - V(M_1)$.
Note that,  for any $j,l \in [k]$,  $|U_1 \cap V_j| = |U_1 \cap V_l|\le k^2$.

If $|U_1 \cap V_1|\leq 1$, then
$\nu(H)\ge |M_1|\ge  n-1$. We may thus assume  $|U_1 \cap V_1|\geq 2$,
and let $S$ be a subset of $U_1$ of type $1$.

Since   $M_1'$ contains at least $4(k+2)\log n \gg k^2$
$S$-absorbing edges, there exists an  $S$-absorbing edge $e\in M_1'$. Thus, $M_2 := (M_1-\{e\}) \cup \{e_S,S_e\}$
is a matching in $H$, with $|M_2|=|M_1|+1$.
Let $U_2 := V(H) - V(M_2)$ and $M_2' := M_1' - \{e\}$.
If $|U_2\cap V_1| \leq 1$, then $\nu(H)\ge |M_2|\ge n-1$. So we may
assume $|U_2\cap V_1| \ge 2$. By repeating this procedure, 
 we get a sequence of matchings $M_1, M_2, \ldots, M_t$ and
a sequence of sets $U_1,U_2, \ldots, U_t$, such that
$|M_{i+1}|=|M_i|+1$, and $U_i=V(H)-V(M_i)$.
Note that $|U_{i+1} \cap V_1|=|U_i\cap V_1|-1$ for $i\in [t-1]$.
Therefore, $t = |U_1 \cap V_1|-1\leq k^2-1$. Hence, $\nu(H)\ge |M_t|\ge n-1$.
\qed


\end{document}